\documentclass[]{amsart}
\usepackage[]{amsmath}
\usepackage{amsthm}
\usepackage{amssymb}
\usepackage{enumitem}

\newtheorem{thm}{Theorem}
\newtheorem{lemma}{Lemma}
\newtheorem{cor}{Corollary}

\theoremstyle{definition}
\newtheorem{prop}{}[section]

\theoremstyle{remark}
\newtheorem*{remark}{Remark}
\newtheorem*{notation}{Notation}

\newcommand{\Osc}[2]{\left(\operatorname{Osc} #1\right)\left( #2 \right)}
\newcommand{\abs}[1]{\left\lvert #1\right\rvert}
\newcommand{\norm}[1]{\left\lVert #1\right\rVert_A}
\newcommand{\Var}{\operatornamewithlimits{Var}}
\newcommand{\card}{\operatorname{card}}
\newcommand{\defas}{\stackrel{\text{\tiny def}}{=}}
\newcommand{\Rp}{\operatorname{Re}}

\title{An Asymptotic Formula for the Sequence  $\norm{\exp(i n h(t))}$}

\author{Bogdan M. Baishanski and Jan Hlavacek}
\address{Bogdan M. Baishanski\\
	 Department of Mathematics\\
	 The Ohio State University\\
	 231 W 18th Avenue\\
	 Columbus,  OH 43210}
\email{bogdan@math.ohio-state.edu}
\address{Jan Hlavacek\\
	 Department of Mathematical Sciences\\
	 Saginaw Valley State University\\
	 7400 Bay Road\\
	 University Center, MI 48710}
\email{jhlavace@svsu.edu}


\begin{document}
\maketitle
\begin{abstract}
Given function $f$ with an absolutely convergent Fourier series, 
\begin{equation*}
   f(t) = \sum_{\nu = -\infty}^{\infty} a_{\nu} e^{i \nu t},
\end{equation*}
define the norm of $f$ as
\begin{equation*}\label{eq_norm}
   \norm{f} = \sum_{\nu = -\infty}^{\infty} \abs{a_\nu}.
\end{equation*}

D.~Girard has in the special case $f(t) = g\left( e^{it} \right)$, $g(z) = \frac{z-\alpha}{1-\overline{\alpha}z}$, $0 < \abs{\alpha} < 1$ proved 
\begin{equation}
   \frac{1}{\sqrt{n}} \norm{f^n} = 
   \frac{8\sqrt{2}}{\Gamma^2\left(\frac{1}{4}\right)} (1-\beta^2)^{\frac{1}{2}}
   F\left(\frac{1}{2},\frac{3}{4};\frac{3}{2};1-\beta^2\right) + o(1), \quad n\rightarrow \infty,
\end{equation}
where $F$ is the hypergeometric function and 
\begin{equation}
   \beta = \frac{1 - \abs{\alpha}}{1 + \abs{\alpha}}.
\end{equation}

In our article we provide a generalization of Girard's formula.

We study the behavior of $\norm{f^n}$ as $n \rightarrow \infty$ 
for $f$ of the form $e^{ih(t)}$,
$h(t + 2\pi) = h(t) + 2k\pi$ for some integer $k$,
and $h$ is a real, odd and  
twice continuously differentiable. We show that if
$h''$ has no zeros in $(0,\pi)$, and 
if it satisfies an additional condition near $0$ and near $\pi$, then
the following asymptotic formula holds:
\begin{equation*}
   \frac{1}{\sqrt{n}}\norm{\exp(i n h)} = 
   \left(\frac{2}{\pi}\right)^{\frac{3}{2}} \int_{0}^{\pi}
   \sqrt{\abs{h''\left(t\right)}}dt
   + o(1) 
\end{equation*}
as $n \rightarrow \infty$.
As one of the corollaries to the result we obtain Girard's formula. Moreover, we show that in the special case when $h$ is $2\pi$-periodic, our formula holds even if the parameter $n$ tends to infinity through real values. 
\end{abstract}
\section{Introduction}

Throughout this article $h$ will denote a real-valued function satisfying $h(t+2\pi) = h(t) + 2k\pi$ for some integer $k$, and such that the Fourier series of 
\[\exp\left(ih(t)\right) = \sum_{\nu=-\infty}^\infty a_\nu e^{i\nu t}\]
is absolutely convergent. With usual notation, 
\[\norm{\exp\left( ih(t) \right)} = \sum_{\nu=-\infty}^\infty \abs{a_\nu}.\]
The basic results on the asymptotic behavior of the sequence $\norm{\exp\left( inh(t) \right)}$ as $n \rightarrow \infty$ were obtained in 1950's. 

The most important of these results is due to A.~Beurling and H.~Helson: from their theorem in \cite{BaH} it follows that
\begin{enumerate}[label=(\arabic*)]
\item If $\norm{\exp\left( inh(t) \right)} = \mathcal{O}(1)$, $n \rightarrow \infty$, then $h(t) = at + b$ where $a$ is an integer and $b$ is a real number. \label{prev:1}
\end{enumerate}
Next in importance are results obtained by J.~P.~Kahane~\cite{kahane_classes}, we quote the following three:
\begin{enumerate}[label=(\arabic*), resume]
\item If $h$ is continuous and piecewise linear, then $\norm{\exp\left( inh(t) \right)} = \mathcal{O}(\log n)$, $n\rightarrow\infty$. (Note that here $\log n$ cannot be replaced by a sequence that would tend to infinity slower, since if $h(x) = \abs{x}$ on $\left(-\pi,\pi\right]$ and $h$ is $2\pi$-periodic, then $\norm{\exp\left( inh(t) \right)} = \frac{2}{\pi}\log n + \mathcal{O}(1)$)\label{prev:2}
\item If on an arbitrarily small interval $h$ is twice differentiable with $h'' \ge \mu > 0$ (or $h'' \le -\mu < 0$), then, for some $\lambda> 0$, $\norm{\exp\left( inh(t) \right)} \ge \lambda\sqrt{n}$.\label{prev:3}
\item If $h$ is analytic and not a linear function, there exist positive constants $\lambda_1$ and $\lambda_2$ such that\label{prev:4}
   \[\lambda_1\sqrt{n} \le \norm{\exp\left( inh(t) \right)} \le \lambda_2\sqrt{n}.\]
\end{enumerate}
Working in a different context B.~Baishanski~\cite{baish} has obtained the following result: 
\begin{enumerate}[label=(\arabic*), resume]
\item If $f$ is analytic in the closed unit disk and $\abs{f\left( e^{it} \right)} = 1$ (i.e. if $f$ is a finite Blaschke product), then \label{prev:5}
   \[\norm{f^n\left( e^{it} \right)} \rightarrow \infty,\quad n\rightarrow \infty,\]
   except in the trivial case $f(z) = cz^m$, $\abs{c} = 1$, $m$ a non-negative integer.
\end{enumerate}
The results we quoted have important applications and \ref{prev:1} has been greatly generalized (for the generalization of \ref{prev:1} and applications of \ref{prev:1} -- \ref{prev:4} see, for example, J.~P.~Kahane \cite{kahane_transforms}, for an application of \ref{prev:5} see P.~Turan \cite{turan}).

It may be not without interest to obtain more precise results on the asymptotic behavior of $\norm{e^{inh}}$.  One such result was obtained by the late Denis Girard. He proved \cite{gir}:
\begin{enumerate}[label=(\arabic*), resume]
\item If 
   \[f(z) = \frac{z-\alpha}{1-\overline{\alpha}z}, \quad 0 < \abs{\alpha} < 1,\]
   then 
   \[
   \frac{1}{\sqrt{n}}\norm{f^n\left( e^{it} \right)} \rightarrow 
   16 \sqrt{2} \left(\Gamma\left(\frac{1}{4}\right)\right)^{-2} \frac{\sqrt{\abs{\alpha}}}{1 +
      \abs{\alpha}}\; F\left(\frac{1}{2}, \frac{3}{4}; \frac{3}{2}; \frac{4\abs{\alpha}}{(1+\abs{\alpha})^2}\right)\]
      as $n\rightarrow\infty$.
\end{enumerate}

\section{Theorems}

We generalize Girard's result by the following 
\begin{thm}
   \label{thm:main}
   Let $h$ be a real, twice continuously differentiable function, $h(t+2\pi) = h(t) + 2k\pi$ for some integer $k$. If, in addition, 
   \begin{enumerate}[label=\textrm{(\arabic*)}]
      \item $h''$ has no zeros on $(0,\pi)$. \label{cond1}
      \item $h$ is odd \label{cond2}
      \item there exist functions $m_0$ and $m_\pi$, positive and increasing on $\left( 0, 4c \right)$ for some $c>0$, and constant $C>0$ such that \label{cond3}
	 \begin{subequations}
	    \begin{equation}
	       \frac{m_0(2t)}{m_0(t)} \le C, \quad \frac{m_\pi(2t)}{m_\pi(t)} \le C, \quad 0 < t < 2c
	       \label{eq:cond_3_part_1}
	    \end{equation}
	    and
	    \begin{gather}
	       m_0(t)  \le h''(t) \le Cm_0(t), \quad 0 < t < 4c, \label{eq:cond_3_part_2}\\
	       m_\pi(t) \le h''(\pi-t) \le Cm_\pi(t), \quad 0 < t < 4c. \label{eq:cond_3_patr_3} 
	    \end{gather}
	    \label{eq:cond_3}
	 \end{subequations}
   \end{enumerate}
   then
   \begin{equation}
      \frac{1}{\sqrt{n}}\norm{e^{inh(t)}} \rightarrow \left( \frac{2}{\pi} \right)^{3/2} \int_0^\pi \sqrt{h''(t)}dt, \quad n\rightarrow \infty
      \label{eq:main}
   \end{equation}
\end{thm}
In the special case when in Theorem~\ref{thm:main} the integer $k$ is equal to $0$, the result \eqref{eq:main} can be strengthened and we have
\begin{thm}
   \label{thm:stronger}
   If $h$ is a $2\pi$-periodic function satisfying the conditions of Theorem~\ref{thm:main}, then
   \begin{equation}
      \frac{1}{\sqrt{x}}\norm{e^{ixh(t)}} \rightarrow \left( \frac{2}{\pi} \right)^{3/2} \int_0^\pi \sqrt{h''(t)}dt,
      \label{eq:main_stronger}
   \end{equation}
   as the real variable $x$ tends to $\infty$.
\end{thm}

\begin{remark}
   If $k\neq 0$ in Theorem~\ref{thm:main}, then $\exp\left( ixh(t+2\pi) \right) \neq \exp\left( ixh(t) \right)$ for every real $t$, unless $x = \frac{m}{k}$ for some integer $m$.  Still we could try to interpret $\norm{\exp(ixh(t))}$ by restricting the function $\exp\left( ixh(t) \right)$ to an interval of length $2\pi$, and then extending that restricted function periodically.  But whatever restriction we choose, the extension will have jump discontinuities.  Therefore its Fourier series will not be absolutely convergent, and therefore for every $x$ such that $xk$ is not an integer, $\norm{\exp(ixh(t))} = \infty$. 
\end{remark}
\begin{remark}
   The conditions \ref{cond1} and \ref{cond2} play an essential role in our proof, since they imply that if $\frac{\nu}{n}$ is in the interior of the range of $h'$, then there are exactly two stationary points in the integral 
   \[a_{n,\nu} = \frac{1}{2\pi} \int_{-\pi}^\pi \exp\left(i(nh(t) - \nu t )\right) dt,\]
   and the contributions from these two points combine nicely. 

   Obviously, $h$ is odd if and only if the coefficients $a_{n,\nu}$ are real. However, we never use directly the fact that $a_{n,\nu}$ are real, neither would our proof be simplified by using the fact directly (indirectly it is reflected in the ``nice'' matching of contributions from the two stationary points).

   Since $h''$ is odd and $2\pi$-periodic, we have that $h''(0) = h''(\pi) = 0$ and the condition \ref{cond3} is a condition on the behavior of $h''$ in a neighbourhood of its zeros. The condition (at $x=0$) is satisfied for example if $h''(t) \sim Ct^\alpha$ or $h''(t) \sim Ct^\alpha\abs{\log t}^\beta$ as $t \rightarrow 0^+$, where $\alpha>0$, $\beta \in \mathbb{R}$. (However, it is not difficult to give examples when the condition \ref{cond3} is not satisfied.  Two such examples are $h''(t) = \exp\left( -\frac{1}{t} \right)$ and $h''(t) = t^2 + t\sin^2\left( \frac{1}{t} \right)$ on some interval $(0,c)$.)
\end{remark}

\begin{notation}
   We denote by $C$ constants, dependent only on the function $h$, not always the same constant.  Since arguments involving $m_0$ and $m_\pi$ are similar, it will be sufficient to consider only $m_0$, which from now on will be denoted $m$. We can assume that $h'' > 0$ on $\left( 0,\pi \right)$, so $h'$ is increasing there from $\alpha = h'(0)$ to $\beta = h'(\pi)$. The inverse function of $h'$ will be denoted by $\psi$.  If $\alpha < \frac{\nu}{n} < \beta$, there exists one and only one $t_{n,\nu}$, $0 < t_{n,\nu} < \pi$ such that $h'(t_{n,\nu}) = 
   \frac{\nu}{n}$ or equivalently $\psi\left( \frac{\nu}{n} \right) = t_{n,\nu}$. 

   We set $\mathcal{M}(x) = \frac{1}{m(\psi(x))}$.  Then $\mathcal{M}$ is defined and decreasing for $\alpha < x < \gamma$, where $\gamma<\beta$ and $\psi(\gamma) < 4c$. The modulus of continuity of $h''$ will be denoted by $\omega$:
   \[\omega(t) = \sup\left\{ \abs{h''(x) - h''(y)}: x,y\in\left[ 0,\pi \right], \abs{x-y} < t \right\}.\]
   Let $\Phi(x)$ denote a function tending to infinity as $x \rightarrow \infty$,
   $\Phi(x) = o(\sqrt{x})$. Additional conditions will be imposed on $\Phi$ in
   the course of the proof.  Finally, $\delta_n = \frac{\Phi(n)}{\sqrt{n}}$,
   and $\alpha_n = \max\left(h'(2\delta_n), \alpha+\frac{1}{n}\right)$, $\beta_n = \min\left(h'(\pi-2\delta_n), \beta-\frac{1}{n}\right)$.
\end{notation}

The proof of Theorem~\ref{thm:main} will be in five steps. 

First we show that 
\[\sum_{\nu\in \mathcal{E}_n} \abs{a_{n,\nu}} = \mathcal{O}(\log n), \quad n \rightarrow \infty,\]
so that the contribution of the external terms $\mathcal{E}_n$, i.e.~the terms $a_{n,\nu}$ such that $\nu \not\in \left( \alpha n+1, \beta n - 1 \right)$ can be ignored in the proof of our formula. 

Second, we establish certain properties of $h''$ needed in the following steps.

Third, we show that 
\[\sum_{\varsigma\in\mathcal{P}_n} \abs{a_{n,\nu}} = \mathcal{O}(1), \quad n \rightarrow \infty\]
where 
\[\mathcal{P}_n = \left\{ \nu: \alpha n + 1 < \nu < \alpha_n n \text{\ or\ } \beta_n n < \nu < \beta n - 1 \right\}\]

In the fourth step we use the method of stationary phase to approximate each of the central terms, that is terms $a_{n,\nu}$ where $\frac{\nu}{n} \in (\alpha_n,\beta_n)$. Finally, in the last step, we use Lemmas~\ref{lem:snk_equiv} and \ref{lem:equidist} to approximate the sum of the central terms and, using this approximation together with estimates obtained in previous steps, we prove the Theorem \ref{thm:main}.

\section{Lemmas}
\begin{lemma}
   Let $\varphi$ be a twice continuously differentiable function on the interval $[a,b]$ such that $\varphi'$ has no zeros in $[a,b]$. The following estimates hold for the integral 
   \[I = \int_a^b e^{i\varphi(t)}dt\]
   \begin{enumerate}[label=\textrm{(\alph*)}]
      \item In the special case $\varphi(b) = \varphi(a) + 2k\pi$ for some integer $k$, and $\varphi'(b) = \varphi'(a)$:
	 \[\abs{I} \le \Var_{[a,b]}\left( \frac{1}{\varphi'} \right)\]\label{lem:lem1_part_a}
      \item In general case, provided $\varphi'$ is monotone on $[a,b]$:
	 \[\abs{I} \le \frac{2}{\min(\abs{\varphi'(a)},\abs{\varphi'(b)})}\]\label{lem:lem1_part_b}
   \end{enumerate}
   \label{lem:lem1}
\end{lemma}
\begin{proof}
   A standard trick is used: the integral is multiplied and divided by $i\varphi'$, then the integration by parts leads to the estimate 
   \[\abs{I} \le \abs{\left.e^{i\varphi(t)}\frac{1}{i\varphi'(t)}\right|_{t=a}^b} + \int_a^b\abs{\frac{\varphi''(t)}{\left( \varphi'(t) \right)^2}},\]
   from which both \ref{lem:lem1_part_a} and \ref{lem:lem1_part_b} follow. 
\end{proof}
\begin{lemma}
   Let $f$ be Riemann integrable on $[c,b]$ for every $c \in (a,b)$, and let the improper integral $\int_{a+}^b f$ exist. 
   Let $a_{j,n} = a + j\frac{b-a}{n}$, $j = 0,\dots, n$, and let $\xi_{j,n}$ be such that $a_{j-1,n} \le \xi_{j,n} \le a_{j,n}$ for $j = 1,\dots,n$.

   If there exists a monotone function $\mathcal{M}$ on $[a,b]$ such that $\abs{f(x)} \le \mathcal{M}(x)$ for $x \in [a,b]$ and $\int_{a+}^b \mathcal{M} < \infty$, then the equidistant Riemann sums of $f$ on $[a,b]$ modified by the omission of the initial summand (i.e.~the sums
   $\frac{b-a}{n}\sum f\left(\xi_{j,n}\right)$ where the summation index does not run from $j=1$ to $n$ but from $j=2$ to $n$)
   converge to $\int_{a+}^b f$.

   In particular, 
   \begin{equation}
      \frac{b-a}{n}\sum f\left( \frac{k}{n} \right) \longrightarrow \int_{a+}^b f
      \label{eq:part_riem_sum}
   \end{equation}
   where the sum runs over all integers $k$ such that $an + 1 < k < bn$.
   \label{lem:part_riem_sum}
\end{lemma}

\begin{remark}
   In the case when $\sum$ denotes the full equidistant Riemann sum, taken over all $k$, $an < k < bn$, \eqref{eq:part_riem_sum} may fail to hold. An example where it fails is provided by the function $f(x) = \frac{1}{\sqrt{x-a}}$ in case $a$ is irrational.
\end{remark}

\begin{proof}
   To show that 
   \begin{equation}
      \mathcal{D}_n = \frac{b-a}{n}\sum_{j=2}^n f\left( \xi_{j,n} \right) - \int_{a+}^b f
      \label{eq:part_riem_proof_1}
   \end{equation}
   tends to zero as $n \rightarrow \infty$, we start, given $\varepsilon > 0$, by choosing $c = c(\varepsilon)$, $d = d(\varepsilon)$, $a < d < c < b$ such that 
   \begin{equation}
      \int_{a+}^c \mathcal{M} < \varepsilon.
      \label{eq:cond_for_c}
   \end{equation}
   Since $f$ is Riemann integrable on $[d,b]$, $\abs{f}$ is bounded there by some constant $K = K(\varepsilon)$.  Let $J = J(n,\varepsilon)$ be defined by the condition $a_{J-1,n} \le c < a_{J,n}$, and let the intervals $I = I(\varepsilon)$ and $I_n = I_n(\varepsilon)$ be defined by $I = [c,b]$ and $I_n = [a_{J,n},b]$. 
   We rewrite \eqref{eq:part_riem_proof_1} as
   \begin{multline*}
      \mathcal{D}_n = \left( \frac{b-a}{n} \sum_{j=J+1}^n f\left( \xi_{j,n} \right) - \int_{a_{J,n}}^b f \right)\\ + \frac{b-a}{n}f\left( \xi_{J,n} \right) - \int_c^{a_{J,n}} f - \int_{a+}^c f + \frac{b-a}{n} \sum_{j=2}^{J-1} f\left( \xi_{j,n} \right)
   \end{multline*}
   so that
   \begin{multline}
      \abs{\mathcal{D}_n} \le \abs{\frac{b-a}{n} \sum_{j=J+1}^n f\left( \xi_{j,n} \right) - \int_{a_{J,n}}^b f }\\ + \frac{b-a}{n}\abs{f\left( \xi_{J,n} \right)} + \int_c^{a_{J,n}}\abs{f} + \int_{a+}^c \abs{f} + \frac{b-a}{n} \sum_{j=2}^{J-1} \abs{f\left( \xi_{j,n} \right)}
      \label{eq:part_riem_proof_est1}
   \end{multline}
   We shall prove the lemma by showing that each of the five summands on the right side of \eqref{eq:part_riem_proof_est1} is less than $\varepsilon$ for $n$ sufficiently large. 

   To estimate the first summand, we shall make use of the following notation: If $J$ is an interval, $\mathcal{P} = \left\{ x_k, 0 \le k \le m \right\}$ a partition of $J$ and $f$ a function on $J$, then 
   \[\Omega\left( f,\mathcal{P}, J \right) = \sum_{k=1}^m \Osc{f}{x_{k-1},x_k} \left( x_k - x_{k-1} \right)\]
   where
   \[\Osc{f}{s,t} = \sup\left\{ \abs{f(x) - f(y)}: s \le x,y \le t \right\}.\]
   We can now write
   \begin{equation}
      \begin{split}
	 \abs{\frac{b-a}{n} \sum_{j=J+1}^n f\left( \xi_{j,n} \right) - \int_{a_{J,n}}^b f } &= \abs{\sum_{j=J+1}^n \int_{a_{j-1,n}}^{a_{j,n}} \left( f\left( \xi_{j,n}\right) - f(x) \right)dx}\\[\baselineskip] &\le \Omega(f,\mathcal{P}_n, I_n),
      \end{split}
      \label{eq:part_riem_proof_est_by_omega}
   \end{equation}
   where $\mathcal{P}_n$ is the partition of the interval $I_n = [a_{J,n},b]$ consisting of points $a_{j,n}$, $J \le j \le n$. 
   Let $\widetilde{\mathcal{P}_n} = \mathcal{P}_n \cup \left\{ c \right\}$.  Then $\widetilde{\mathcal{P}_n}$ is a partition of $I = [c,b]$ and $\left\lVert \widetilde{\mathcal{P}_n}\right\rVert = \frac{b-a}{n}$.  Obviously
   \begin{equation}
      \Omega\left( f,\mathcal{P}_n, I_n \right) \le \Omega\left( f, \widetilde{\mathcal{P}_n},I\right).
      \label{eq:extending_omega}
   \end{equation}
   Since $f$ is Riemann integrable on $I$, $\Omega\left( f, \widetilde{\mathcal{P}_n},I\right) \rightarrow 0$ as $\left\lVert \widetilde{\mathcal{P}_n}\right\rVert \rightarrow 0$. It follows from \eqref{eq:part_riem_proof_est_by_omega} and \eqref{eq:extending_omega} that the first summand in \eqref{eq:part_riem_proof_est1} is less than $\varepsilon$ for large $n$. 

   The second summand is less than $\frac{b-a}{n}K(\varepsilon)$, the same holds for the third summand. The fourth summand is less then $\varepsilon$ by \eqref{eq:cond_for_c}.  The crucial step in the proof is the estimate of the fifth summand. 

   Since $\mathcal{M}$ is monotone decreasing on $(a,b]$ we have, on one hand
   \[\sum_{j=2}^{J-1} \abs{f\left(\xi_{j,n}\right)} \le \sum_{j=2}^{J-1} \mathcal{M}\left(\xi_{j,n}\right) \le \sum_{j=2}^{J-1} \mathcal{M}\left(a_{j-1,n}\right) = \sum_{j=1}^{J-2} \mathcal{M}\left(a_{j,n}\right) \]
   and on the other hand 
   \[\frac{b-a}{n}\mathcal{M}\left( a_{j,n} \right) \le \int_{a_{j-1,n}}^{a_{j,n}} \mathcal{M}(x) dx\]
   and so
   \[\frac{b-a}{n} \sum_{j=2}^{J-1} \abs{f\left( \xi_{j,n} \right)} \le \int_{a+}^{a_{J-2,n}} \mathcal{M}(x) dx \le \int_{a+}^c \mathcal{M}(x) dx \le \varepsilon\]
   by \eqref{eq:cond_for_c}.
\end{proof}

\begin{lemma}
   Given an interval $J_0$ and an array $s_{n,k}$, $0 \le s_{n,k} \le 1$, where
   $n = 1, 2, \dots$, and $k$ runs through integers such that $\frac{k}{n} \in
   J_0$, the following three conditions are equivalent:
   \begin{enumerate}[label=\textrm{(\alph*)}]
      \item For any pair of intervals $I$, $J$, $I \subset [0,1]$, $J \subset J_0$, we have, as $n \rightarrow \infty$,
	 \[\frac{\mathcal{N}_n}{n} \rightarrow \abs{I}\cdot\abs{J},\]
	 where $\mathcal{N}_n$ is the number of indices $k$ such that $s_{n,k} \in I$ and $\frac{k}{n} \in J$, and $\abs{I}$ and $\abs{J}$ denote the lengths of the intervals $I$ and $J$, respectively.\label{cond:l3a}
      \item For any non-zero integer $j$ and any subinterval $J$ of $J_0$\label{cond:l3b}
	 \[\frac{1}{n} \sum_{\left\{k:\frac{k}{n}\in J  \right\}} \exp\left( 2\pi ijs_{n,k} \right) \rightarrow 0, n \rightarrow \infty\]
      \item If $f$ is Riemann integrable on the interval $J_0$, and $g$ is Riemann integrable on $[0,1]$, then\label{cond:l3c}
	 \[\frac{1}{n} \sum_{\left\{ k:\frac{k}{n} \in J_0 \right\}} f\left( \frac{k}{n} \right) g(s_{n,k}) \rightarrow \int_{J_0} f \cdot \int_0^1 g, \quad n\rightarrow \infty\]
   \end{enumerate}
   \label{lem:snk_equiv}
\end{lemma}
\begin{proof}
   The proof is an adaptation of the well known method of Hermann Weyl.
   Let 
   \begin{align*}
   (f,g)_n & \defas \frac{1}{n}\sum_{k: \frac{k}{n}\in J_0} f\left( \frac{k}{n} \right)g(s_{n,k})\\
   (f,g) & \defas \int_{J_0} f \cdot \int_0^1 g
   \end{align*}
   and let $\mathcal{F}$ denotes the collection of pairs $\left\{ f,g \right\}$ such that $(f,g)_n \rightarrow (f,g)$ as $n\rightarrow \infty$.  We observe that 
   \begin{itemize}
      \item[(*)] If $f$ is Riemann integrable on $J_0$ then $\left\{ f,1 \right\} \in \mathcal{F}$. \label{fact:riem}
   \end{itemize}

   It is easy to verify that condition \ref{cond:l3a} can be phrased as
   \begin{enumerate}[label=(\alph*')]
      \item \label{cond:l3aprime} $\left\{\chi_J, \chi_I \right\} \in \mathcal{F}$ for any pair of intevals $I$, $J$, $I \subset [0,1]$, $J \subset J_0$
   \end{enumerate}
   and conditions \ref{cond:l3b} and \ref{cond:l3c} as
   \begin{enumerate}[label=(\alph*'), resume]
      \item \label{cond:l3bprime} $\left\{ \chi_J, g_j \right\} \in \mathcal{F}$ for $g_j(x) = \exp\left( 2\pi ijx \right)$, $j=\pm1, \pm2, \dots$, and $J$ any subinterval of $J_0$,
      \item \label{cond:l3cprime} $\left\{ f,g \right\} \in \mathcal{F}$ for any complex-valued functions $f$ and $g$, that are Riemann-integrable on $J_0$ or $[0,1]$, respectively. 
   \end{enumerate}
   It is obvious that \ref{cond:l3aprime} and \ref{cond:l3bprime} are consequences of \ref{cond:l3cprime}, so to prove the equivalence of the three conditions it will be sufficient to show that each of \ref{cond:l3aprime} and \ref{cond:l3bprime} implies \ref{cond:l3cprime}.  This will be done in two steps --- first we show that \ref{cond:l3bprime} implies a condition slightly weaker than \ref{cond:l3aprime}, which we denote by \ref{cond:l3adp}, and then we show that \ref{cond:l3adp} implies \ref{cond:l3cprime}:
   \begin{enumerate}[label=(\alph*'')]
      \item \label{cond:l3adp} $\left\{ \chi_J, \chi_I \right\} \in \mathcal{F}$ for any pair of intervals $I$, $J$ such that $J \subset J_0$ and that the closure of the interval $I$ lies in the open interval $(0,1)$. 
   \end{enumerate}
   The proof will use the following facts, which are easily established:
   \begin{itemize}
      \item[(**)] If $\left\{ f,\varphi_k \right\}\in \mathcal{F}$ and $\varphi_k$ converge to $\varphi$ uniformly on $[0,1]$, then $\left\{ f,\varphi \right\} \in \mathcal{F}$. \label{fact:unif}
      \item[(***)'] If $f \ge 0$ on $J_0$ and $g$ is real-valued and if, for every $\varepsilon>0$, there exist functions $\varphi_1$ and $\varphi_2$ on $[0,1]$, $\varphi_1 \le g \le \varphi_2$,  such that $\left\{ f, \varphi_1 \right\} \in \mathcal{F}$, $\left\{ f,\varphi_2 \right\} \in \mathcal{F}$ and $\int_0^1(\varphi_2 - \varphi_1) < \varepsilon$, then $\left\{ f,g \right\}\in \mathcal{F}$.\label{fact:approxg}
      \item[(***)''] Similarly, if $g\ge0$ and $f$ is real-valued, and if for
	 every $\varepsilon>0$ there exist functions $\varphi_1$ and
	 $\varphi_2$ on $J_0$ such that $\varphi_1 \le f \le \varphi_2$,
	 $\left\{ \varphi_1, g \right\} \in \mathcal{F}$, $\left\{ \varphi_2, g
	 \right\}\in \mathcal{F}$ and $\int_{J_0}(\varphi_2 - \varphi_1) <
	 \varepsilon$, then $\left\{ f,g \right\} \in \mathcal{F}$.\label{fact:approxf}
   \end{itemize} 

   We start the proof by observing that (*) and \ref{cond:l3bprime} imply that
   $\left\{ \chi_j, T \right\} \in \mathcal{F}$ for any trigonometric
   polynomial $T(x) = \sum_{\abs{\nu}\le N} c_\nu e^{2\pi i\nu x}$.  From (**)
   we obtain that $\left\{ \chi_J, \varphi \right\} \in \mathcal{F}$ for any
   continuous function $\varphi$ on $[0,1]$ satisfying $\varphi(0) =
   \varphi(1)$. Since for any $\varepsilon>0$ there exist such continuous
   functions $\varphi_1$, $\varphi_2$ with the property $\varphi_1 \le \chi_I
   \le \varphi_2$ and $\int_0^1 (\varphi_2 - \varphi_1) < \varepsilon$ (since
   the endpoints of $I$ are inside $(0,1)$), we obtain by (***)' that $\left\{ \chi_j,\chi_I  \right\} \in \mathcal{F}$, i.e~\ref{cond:l3adp} is proved. 

   From (*) we get $\left\{ \chi_J, 1 \right\} \in \mathcal{F}$.  Together with \ref{cond:l3adp} this implies that $\left\{ \chi_J, \varphi \right\} \in \mathcal{F}$ for any step-function $\varphi$ satisfying $\varphi(0) = \varphi(1)$.  If $g$ is a real-valued Riemann-integrable function on $[0,1]$, for each $\varepsilon>0$ we can find such step-functions $\varphi_1$, $\varphi_2$ that $\varphi_1 \le g \le \varphi_2$ and $\int_0^1 (\varphi_2 - \varphi_1) < \varepsilon$. Therefore, using (***)' again, we obtain $\left\{ \chi_J, g \right\} \in \mathcal{F}$, and $\left\{ \varphi, g \right\} \in \mathcal{F}$ for $\varphi$ a real-valued step-function on $J_0$.  In particular this holds for any $g$ non-negative, Riemann-integrable.  Using now (***)'' we obtain that $\left\{ f,g \right\} \in \mathcal{F}$ for any $f$ and $g$ real-valued, Riemann-integrable, $g$ non-negative.  By linearity, the same holds if both $f$ and $g$ are complex-valued. This gives us \ref{cond:l3cprime}, which concludes the proof of the lemma.
\end{proof}

\begin{remark}
   If an array $s_{n.k}$ satisfies any of the three conditions in Lemma~\ref{lem:snk_equiv} (and therefore all of these conditions), we say that the array is equidistributed on the interval $[0,1]$.
\end{remark}

\begin{lemma}
   Let $\varphi$ be a real-valued twice differentiable function on an interval $J_0$ and let $\varphi''(x) \ge \rho > 0$ for $x \in J_0$.  Then the array $\left\langle n\varphi\left( \frac{k}{n} \right)\right\rangle$, where $\left\langle s \right\rangle$ denotes the fractional part of $s$, for $n = 1, 2, \dots$ and $\frac{k}{n} \in J_0$, is equidistributed on $[0,1]$.
   \label{lem:equidist}
\end{lemma}
\begin{proof}
   It is suficient to show that the condition \ref{cond:l3b} of Lemma~\ref{lem:snk_equiv} is satisfied. Let $j$ denote any integer different from zero. Since $\exp\left( 2\pi ijx \right)$ has period 1, we get for any subinterval $J = [c,d]$ of $J_0$
   \begin{equation}
      \mathcal{U}_{n,j} \defas \sum_{k|\frac{k}{n}\in J} \exp\left( 2\pi ijs_{n,k} \right) = \sum_{nc \le k \le nd} \exp\left(2\pi ijn\varphi\left( \frac{k}{n} \right)\right)
      \label{eq:unjdef}
   \end{equation}
   To estimate the last expression, we use one of Van der Corput's lemmas \cite[ch.~5, lemma 4.6]{zygmund1}:
   \begin{itemize}
      \item[(*)] If $f''(x) \ge \mu > 0$ of $f''(x) \le -\mu < 0$ on $[a,b]$, then \label{lem:van_der_corput}
	 \begin{equation}
	    \abs{\sum_{a < k \le b} \exp(2\pi if(k))} \le \left( \abs{f'(b) - f'(a)} + 2 \right)\left( \frac{4}{\sqrt{\mu}} + A \right)
	    \label{eq:van_der_corput}
	 \end{equation}
	 where $A$ is an absolute constant. 
   \end{itemize}
   We shall apply this lemma to the case $f(x) = jn\varphi\left( \frac{x}{n} \right)$, $a = nc$, $b = nd$.  Since $f'(x) = j\varphi'\left( \frac{x}{n} \right)$, $\abs{f''(x)} = \frac{\abs{j}}{n}\varphi''\left( \frac{x}{n} \right) \ge \frac{\abs{j}}{n}\rho$, the conditions of (*) are satisfied with $\mu = \frac{\abs{j}}{n}\rho$. Therefore we obtain from \eqref{eq:van_der_corput} that 
   \[
   \abs{\sum_{nc < k \le nd} \exp\left( 2\pi ijn\varphi\left( \frac{k}{n} \right) \right)} \le \abs{j\varphi'(d) - j\varphi'(c)}\left( \frac{4}{\sqrt{\frac{\abs{j}}{n}\rho}} + A \right) \le C_j\sqrt{n} + D_j,
   \]
   where $C_j$ and $D_j$ are constant depending only on $j$ and the function $\varphi$. 

   Then, by \eqref{eq:unjdef},
   \[
   \frac{\abs{\mathcal{U}_{n,j}}}{n} \le \frac{C_j}{\sqrt{n}} + \frac{D_j}{n} \rightarrow 0
   \]
   as $n \rightarrow \infty$, and therefore the condition \ref{cond:l3b} of Lemma~\ref{lem:snk_equiv} is satisfied.
\end{proof}
\section{External Terms}

In this section we will provide an estimate for the sum of the external terms, namely
\begin{equation}
   \sum_{\nu\in \mathcal{E}_n} \abs{a_{n,\nu}} = \mathcal{O}(\log n), \quad n \rightarrow \infty
   \label{eq:sum_external}
\end{equation}
where the external terms $\mathcal{E}_n$ are the terms where $\frac{\nu}{n} \le \alpha + \frac{1}{n}$ or $\frac{\nu}{n} > \beta-\frac{1}{n}$.  We will provide an estimate for the sum of the terms where $\frac{\nu}{n} \le \alpha + \frac{1}{n}$, the estimate for the other case can be obtained in a similar way. 

We start by splitting the sum into three parts:
\begin{equation}
   \sum_{\frac{\nu}{n} \le \alpha + \frac{1}{n}} \abs{a_{n,\nu}} = 
   \sum_{\alpha - \frac{2}{n} < \frac{\nu}{n} \le \alpha + \frac{1}{n}} \abs{a_{n,\nu}}  +
   \sum_{\alpha - 1 < \frac{\nu}{n} \le \alpha - \frac{2}{n}} \abs{a_{n,\nu}}  +
   \sum_{\frac{\nu}{n} \le \alpha - 1} \abs{a_{n,\nu}}  
   \label{eq:external_split}
\end{equation}
We can easily dismiss the first sum, as it has at most three terms, each less than or equal to one. 

For the second sum, we will use the following estimate of $\abs{a_{n,\nu}}$:
\begin{equation}
   \abs{a_{n,\nu}} = \frac{1}{2\pi}\abs{b_{n,\nu} + \overline{b_{n,\nu}}} \le \frac{1}{\pi}\abs{b_{n,\nu}}
   \label{eq:estim_abs_annu_with_bnnu}
\end{equation}
where
\begin{equation}
   b_{n,\nu} = \int_0^\pi e^{i(nh(t) - \nu t)} dt.
   \label{eq:bnnu_def}
\end{equation}
Applying Lemma~\ref{lem:lem1} part \ref{lem:lem1_part_b} with $\varphi(t) = nh(t) - \nu t$ we obtain 
\begin{equation}
   \abs{b_{n,\nu}} \le \frac{2}{\min\left( \abs{n\alpha-\nu},\abs{n\beta-\nu} \right)}
   \label{eq:bnnu_estim}
\end{equation}
This gives us estimate
\begin{equation}
   \sum_{\alpha-1 \le \frac{\nu}{n} \le \alpha-\frac{2}{n}} \abs{b_{n,\nu}} \le 
   \sum_{\alpha-1 \le \frac{\nu}{n} \le \alpha-\frac{2}{n}} \frac{2}{\alpha n - \nu} \le C\log(n)
   \label{eq:bnnu_sum_estimate}
\end{equation}

For the third summand we apply Lemma~\ref{lem:lem1} part \ref{lem:lem1_part_a}, again with $\varphi(t) = nh(t) - \nu t$, but with the interval $(a,b) = (-\pi,\pi)$.  It is easy to verify that if $\frac{\nu}{n} \le \alpha-1$, $\varphi$ satisfies the hypothesis of Lemma~\ref{lem:lem1} part \ref{lem:lem1_part_a}, and we obtain
\begin{equation}
   \abs{a_{n,\nu}} = \frac{1}{2\pi}\abs{\int_{-\pi}^\pi e^{i(nh(t) - \nu t)} dt} \le \frac{1}{2\pi}\Var_{(-\pi,\pi)}\left( \frac{1}{\varphi'} \right).
   \label{eq:annu_estim_var}
\end{equation}
Since $\varphi'$ is even and monotone on $(0,\pi)$, we have
\[\Var_{(-\pi,\pi)}\left( \frac{1}{\varphi'} \right) = 2\left(
\frac{1}{\varphi'(0)} - \frac{1}{\varphi'(\pi)} \right) \le \frac{2}{n}\frac{\beta-\alpha}{\left( \alpha-\frac{\nu}{n} \right)^2}\]
which together with \eqref{eq:annu_estim_var} gives
\begin{equation}
   \sum_{\frac{\nu}{n} < \alpha-1}\abs{a_{n,\nu}} \le C\sum_{\frac{\nu}{n} < \alpha-1}\frac{1}{n}\frac{\beta-\alpha}{\left( \alpha-\frac{\nu}{n} \right)^2} \le C\int_{-\infty}^{\alpha-1} \frac{\beta-\alpha}{(\alpha-x)^2}dx
   \label{eq:estim_far_terms}
\end{equation}
which is a constant independent of $n$.

Then \eqref{eq:external_split} together with \eqref{eq:estim_abs_annu_with_bnnu}, \eqref{eq:bnnu_sum_estimate} and \eqref{eq:estim_far_terms} will give us \eqref{eq:sum_external}. 

\section{Some Auxiliary Results on $h''$ and $\sqrt{h''}$}
\begin{prop}\label{prop:first}
The function $\frac{1}{h''\left( \psi(t) \right)}$ is integrable
on $\left[
\alpha,\beta \right] = \left[ h'(0), h'(\pi) \right]$ (and therefore
$\frac{1}{\sqrt{h''\left( \psi(t) \right)}}$ is also integrable on $\left[
\alpha,\beta \right]$). 
\end{prop}
   To prove this fact, we consider the integral 
\[
I(\tau_1, \tau_2) = \int_{h'(\tau_1)}^{h'(\tau_2)} \frac{1}{h''\left(
\psi(t) \right)} dt, \quad 0 < \tau_1 < \tau_2 < \pi.
\]
Making the change of variable $u = \psi(t)$ we get $t = h'(u)$, $dt = h''(u)
du$ and so $I(\tau_1,\tau_2) = \int_{\tau_1}^{\tau_2} du = \tau_2 - \tau_1$,
which converges as $\tau_1 \rightarrow 0+$, $\tau_2 \rightarrow \pi-$. 

For further reference, we note
\begin{equation}
   \label{eq:int_of_1_over_hdblprime}
   \int_{h'(0)}^{h'(d)} \frac{1}{h''(\psi(t))} dt = d, \quad 0 \le d \le \pi
\end{equation}
If $d>0$ is chosen sufficiently small so that $h''(x) \le 1$ for $0 < x \le d$,
then we deduce from~\eqref{eq:int_of_1_over_hdblprime} that 
\begin{equation}
   \int_{h'(0)}^{h'(d)} \frac{1}{\sqrt{h''(\psi(t))}} dt \le d.
   \label{eq:int_1_over_sqrt_estim}
\end{equation}

Finally, we note that the same change of variable as above ($u = \psi(t)$)
brings another simplification:
\begin{equation}
   \int_{h'(\tau_1)}^{h'(\tau_2)} \frac{1}{\sqrt{h''\left( \psi(t)
   \right)}}dt = \int_{\tau_1}^{\tau_2} \sqrt{h''(u)} du.
   \label{eq:int_1_over_sqrt_equals}
\end{equation}

\begin{prop}\label{prop:second}
   The function $\frac{1}{h''(\psi(t))}$ has on $\left[ h'(0), h'(4c)
   \right]$ a monotone majorant $\mathcal{M}(t)$ which is integrable.
   (Therefore $\frac{1}{\sqrt{h''(\psi(t))}}$ has also such a majorant:
   $\sqrt{\mathcal{M}(t)}$). (Both statements are valid also on the interval
   $\left[ h'(\pi-4c), h'(\pi) \right]$, possibly with different majorant
   $\mathcal{M}(t)$).
\end{prop}
Clearly, since $h''(x) \ge m(x)$ on $\left[ 0,4c \right]$, we have that
$\mathcal{M}(t) = \frac{1}{m(\psi(t))}$ is a monotone majorant of
$\frac{1}{h''(\psi(t))}$.  The function $\mathcal{M}$ is integrable because
it is itself majorized by the integrable function $\frac{C}{h''(\psi(t))}$:
this follows from the inequality $h''(x) \le Cm(x)$, $x \in \left[ 0,4c
\right]$. 

\begin{prop}
   According to the paragraphs~\ref{prop:first} and~\ref{prop:second}, conditions
   of Lemma~\ref{lem:part_riem_sum} are satisfied by the functions
   $\frac{1}{h''\left( \psi(t) \right)}$ and $\frac{1}{\sqrt{h''\left(
   \psi(x) \right)}}$ on the interval $\left[ 0,4c \right]$.  By
   Lemma~\ref{lem:part_riem_sum} we have then
   \begin{equation}
      \sum_{nh'(0) + 1 \le k \le nh'(d)}\frac{1}{n} \frac{1}{h''\left( \psi\left(
      \frac{k}{n} \right) \right)} \le C \int_0^d \mathcal{M}(x) dx, 0 < d\le
      h'(4c).
      \label{eq:third_one}
   \end{equation}
   Since $\mathcal{M}(x) \le \frac{C}{h''\left( \psi(x) \right)}$ for $x \in
   \left( h'(0), h'(4c) \right]$, we deduce
   from~\eqref{eq:int_of_1_over_hdblprime} and~\eqref{eq:third_one} that
   \begin{equation}
      0 \le \sum_{nh'(0) + 1 \le k \le nh'(d)} \frac{1}{n}\frac{1}{h''\left( \psi\left(
      \frac{k}{n} \right) \right)} \le Cd
      \label{eq:par_3_estim1}
   \end{equation}
   if $0 < d \le 4c$. 

   Similarly, using~\eqref{eq:int_1_over_sqrt_estim} we obtain
   \begin{equation}
      0 \le \sum_{nh'(0) + 1 \le k \le nh'(d)} \frac{1}{n}\frac{1}{\sqrt{h''\left( \psi\left(
      \frac{k}{n} \right) \right)}} \le Cd
      \label{eq:par_3_estim1_for_sqrt}
   \end{equation}
   if $0 \le d \le \min\left( 4c, d_0 \right)$ where $h''(x) \le 1$ for
   $0 < x \le d_0$.
\end{prop}
\begin{prop}
   In this paragraph we prove the following inequality to which we shall refer
   repeatedly:
   \begin{equation}
      h''(t) \le Ch''(\xi)
      \label{eq:useful_hdoubleprime_inequality}
   \end{equation}
   in each of the following cases:
   \begin{enumerate}[label=\textrm{(\alph*)}]
      \item $0 < t < \xi < 4c$ \label{case:useful_ineq_c1}
      \item $\pi-4c < \xi < t < \pi$ \label{case:useful_ineq_c2}
      \item if $\delta$ is any positive number less than $c$ and if $2\delta <
	 t < \pi-2\delta$, $\abs{\xi-t} < \delta$.  \label{case:useful_ineq_c3}
   \end{enumerate}
   Note that the constant $C$
   is~\eqref{eq:useful_hdoubleprime_inequality} does not depend on
   $\delta$.  
\end{prop}
\noindent To prove the inequality, we will consider each of the three cases.
   \begin{enumerate}[label=\textrm{(\alph*)}]
      \item \label{case:a} By monotonicity of $m$ and the inequality $m \le h'' \le Cm$ on
	 $\left[ 0,4c \right]$ we get $h''(t) \le Cm(t) \le Cm(\xi) \le
	 Ch''(\xi)$.
      \item is analogous to \ref{case:a}.
      \item This case we need to subdivide into three subcases:
	 \begin{enumerate}[label=\textrm{(\alph{enumi}\arabic*)}]
	    \item \label{case:c1} $2c \le t \le \pi-2c$, $\abs{\xi-t}\le \delta < c$
	    \item \label{case:c2} $2\delta \le t \le 2c$, $\abs{\xi-t}\le \delta < c$
	    \item \label{case:c3} $\pi-2c \le t \le \pi-2\delta$, $\abs{\xi-t}\le \delta < c$
	 \end{enumerate}
	 The case \ref{case:c1} is simple: If $I = \left[ 2c, \pi-2c \right]$
	 and $J = \left[ c, \pi-c \right]$, it is sufficient to choose $C$ so
	 that $\sup_I h'' < C \inf_J h''$. 

	 It remains only to consider the case \ref{case:c2}, because
	 \ref{case:c3} can be treated the same way.  This is the only delicate
	 case, and this is the only place in the proof of the Theorem \ref{thm:main} where we
	 use the assumption \eqref{eq:cond_3_part_1}, that is
	 $\frac{m(2t)}{m(t)} \le C$ for $0 \le t \le 2c$. 

	 We observe that $0 \le t-\delta < \xi < t + \delta < 3c$, and so 
	 \begin{equation}
	    m(t-\delta) \le m(\xi) \le h''(\xi).
	    \label{eq:Y}
	 \end{equation}
	 Since $t-\delta < t < 2c$, the assumption \eqref{eq:cond_3_part_1} gives
	 us $m(2(t-\delta)) \le Cm(t-\delta)$ and from \eqref{eq:Y} it follows
	 that $m(2(t-\delta)) \le Ch''(\xi)$. On the other hand, $t \ge 2\delta$
	 implies $t \le 2(t-\delta)$, and so $m(t) \le m(2(t-\delta))$. 
	 Now $h''(t) \le Cm(t) \le Cm(2(t-\delta)) \le Ch''(\xi)$, so
	 \eqref{eq:useful_hdoubleprime_inequality} holds. 
   \end{enumerate}

\section{Total contribution of periphery terms}
We will now prove that 
\begin{equation}
   \sum_{\nu\in \mathcal{P}_n} \abs{a_{n,\nu}} = \mathcal{O}(1),\quad n \rightarrow \infty.
   \label{eq:periphery_sum}
\end{equation}
Here 
\[
\begin{aligned}
   \mathcal{P}_n &= \mathcal{P}_n^l \cup \mathcal{P}_n^r,\\
   \mathcal{P}_n^l &= \left\{ \nu: \alpha n + 1 < \nu < \alpha_n n \right\},\\
   \mathcal{P}_n^r &= \left\{ \nu: \beta_n n < \nu < \beta n - 1\right\},\\
   \alpha_n &= \max\left( h'(2\delta_n), h'(0)+\frac{1}{n} \right),\\
   \beta_n &= \min\left( h'(\pi-2\delta_n), h'(\pi) - \frac{1}{n} \right).
\end{aligned}
\]
(We should mention that one of the sets $\mathcal{P}_n^l$ or $\mathcal{P}_n^r$, or both of them, may be empty).
We shall consider only the left peripheral terms, the right periphery terms can be treated similarly. 

It will be sufficient to consider the integrals
\[ b_{n,\nu} = \int_0^\pi \exp\left( i(nh(t) - \nu t) \right) dt,\]
where $\alpha n + 1 < \nu < \alpha_n n$, or equivalently $h'(0) + \frac{1}{n} < \frac{\nu}{n} < h'(2\delta_n)$.  We write
\[ b_{n,\nu} = \int_0^{3\delta_n} \exp\left( i(nh(t) - \nu t) \right) dt + \int_{3\delta_n}^\pi \exp\left( i(nh(t) - \nu t) \right) dt .\]
For the first summand we use the trivial estimate
\begin{equation}
   \abs{\int_0^{3\delta_n} \exp\left( i(nh(t) - \nu t) \right) dt} \le 3\delta_n.
   \label{eq:periph_first_summand}
\end{equation}
To estimate the second summand, we shall apply Lemma~\ref{lem:lem1}, part~\ref{lem:lem1_part_b}.  We set $\varphi(t) = nh(t) - \nu t$ and observe that since $h'' > 0$ on $[0,\pi]$, the function $\varphi'$ is strictly increasing.  Observing that $\nu \in \mathcal{P}_n^l$ implies $h'\left(t_{n,\nu}\right) = \frac{\nu}{n} < h'(2\delta_n)$, we have $t_{n,\nu} < 2\delta_n$. Since $\varphi'(t_{n,\nu}) = 0$, we get that $\varphi'(t) \ge \varphi'(3\delta_n) > 0$ on $[3\delta_n, \pi]$.  So, by Lemma~\ref{lem:lem1}, part~\ref{lem:lem1_part_b}, we get that
\begin{equation}
   \abs{\int_{3\delta_n}^\pi \exp\left( i(nh(t) - \nu t) \right) dt } \le \frac{2}{\varphi'(3\delta_n)} = \frac{2}{n\left( h'(3\delta_n) - \frac{\nu}{n} \right)}. 
   \label{eq:periph_second_summand_1}
\end{equation}
But $\frac{\nu}{n} = h'(t_{n,\nu})$ and
\begin{equation}
   h'(3\delta_n) - \frac{\nu}{n} = h'(3\delta_n) - h'\left(t_{n.\nu}\right) = h''(\xi)\left( 3\delta_n - t_{n,\nu} \right)
   \label{eq:periph_second_summand_2}
\end{equation}
for some $\xi \in (t_{n,\nu},3\delta_n)$.  Since $t_{n,\nu} < 2\delta_n$, we have $3\delta_n - t_{n,\nu} \ge \delta_n$.  
On the other hand, by the inequality \eqref{eq:useful_hdoubleprime_inequality}, case \ref{case:useful_ineq_c1}, since $0 < t_{n,\nu} < \xi < 3\delta_n < 4c$, we have $h''(\xi) \ge \frac{1}{C} h''\left( t_{n,\nu} \right)$. Thus we derive from \eqref{eq:periph_second_summand_2} that
\[ h'(3\delta_n) - \frac{\nu}{n} \ge \frac{\delta_n h''\left(t_{n,\nu}\right)}{C}.\]
From this estimate and \eqref{eq:periph_second_summand_1} we obtain
\[\abs{\int_{3\delta_n}^\pi \exp\left( i(nh(t) - \nu t) \right) dt } \le \frac{C}{n\delta_nh''\left( t_{n,\nu} \right)}.\]
The last estimate and \eqref{eq:periph_first_summand} imply
\begin{equation}
   \sum_{\nu\in \mathcal{P}_n^l} \abs{b_{n,\nu}} \le 3\delta_n \sum_{\nu\in\mathcal{P}_n^l} 1 + \frac{C}{\delta_n} \sum_{\nu\in\mathcal{P}_n^l} \frac{1}{nh''\left( t_{n,\nu} \right)}.
   \label{eq:periph_both_summands}
\end{equation}
We have 
\[
\begin{aligned}
\delta_n \sum_{\nu\in\mathcal{P}_n^l} 1 &= \delta_n\cdot\card\left\{ \nu: \alpha n + 1 < \nu < \alpha_n n\right\}\\ 
&\le \delta_n\cdot \left( \alpha_n - \alpha \right)\cdot n\\ 
&= \left( h'(2\delta_n) - h'(0) \right)\cdot n \delta_n\\ 
&= h''(\xi)\cdot 2n\delta_n^2 
\end{aligned}
\]
for some $\xi\in(0,2\delta_n)$. Since $h''(x) \le Cm(x)$ and $m$ is increasing, we get
\[ h''(\xi)\cdot 2n\delta_n^2 \le C m(\xi)n\delta_n^2 \le C m(2\delta_n)n\delta_n^2 \]
Since, on the other hand, $m(x) \le h''(x) = h''(x) - h''(0) \le \omega(x)$, and $\omega(2x) \le 2 \omega(x)$, we have
\[C m(2\delta_n)n\delta_n^2 \le C\omega(2\delta_n)n\delta_n^2 \le C\omega(\delta_n)n\delta_n^2\]
This, together with $\delta_n = \frac{\Phi(n)}{\sqrt{n}}$, and our choice of
$\Phi(n)$
which will be made in \eqref{eq:choice_of_Phi_n}, namely $\omega(\delta_n) \Phi(n)^4 = 1$, gives us
\begin{equation}
 \delta_n \sum_{\nu\in\mathcal{P}_n^l} 1 \le C\omega(\delta_n)\Phi(n)^2 = \frac{C}{\Phi(n)^2} = o(1),\quad n\rightarrow \infty.
   \label{eq:periph_first_summamnd_again}
\end{equation}

We will now show that the second term on the right hand side of \eqref{eq:periph_both_summands} is bounded. 
Using \eqref{eq:cond_3_part_2} we get
\[\sum_{\nu \in \mathcal{P}_n^l} \frac{1}{nh''\left( \psi\left( \frac{\nu}{n} \right) \right)} \le \frac{1}{n}\sum_{\alpha n + 1 < \nu < \alpha_n n} \frac{1}{m\left(\psi\left( \frac{\nu}{n} \right)\right)}\]
That can then be written as
\[\frac{1}{n} \sum \mathcal{M}\left( \frac{\nu}{n} \right)\] 
where the sum runs over $\nu$ satisfying $h'(0) + \frac{1}{n} \le \frac{\nu}{n} \le h'(2\delta_n)$.
Since $\mathcal{M}$ is decreasing, this sum is less than or equal to
\[
\begin{aligned}
\int_{h'(0)}^{h'(2\delta_n)} \mathcal{M}(x) dx &= \int_{h'(0)}^{h'(2\delta_n)} \frac{1}{m(\psi(x))} dx\\ 
&\le C \int_{h'(0)}^{h'(2\delta_n)} \frac{1}{h''(\psi(x))} dx = 2\delta_n C
\end{aligned}
\]
where the last step follows from \eqref{eq:int_of_1_over_hdblprime}.
Thus
\begin{equation}
   \frac{1}{\delta_n}\sum_{\nu\in \mathcal{P}_n^l} \frac{1}{nh''\left( t_{n,\nu} \right)} \le C. 
   \label{eq:periph_second_summand_again}
\end{equation}

From the previous estimates \eqref{eq:periph_both_summands}, \eqref{eq:periph_first_summamnd_again} and \eqref{eq:periph_second_summand_again}, and the observation \eqref{eq:estim_abs_annu_with_bnnu} that $\abs{a_{n,\nu}} \le \abs{b_{n,\nu}}$ we get 
\[\sum_{\nu\in \mathcal{P}_n^l} \abs{a_{n,\nu}} = \mathcal{O}(1), \quad n\rightarrow \infty.\]
Since similar estimate holds for $\mathcal{P}_n^r$, we obtain \eqref{eq:periphery_sum}.

\section{Central Terms}
In this section we show that, for $\nu\in \mathcal{C}_n$, i.e.~$n\alpha_n \le
\nu \le n\beta_n$, 
\begin{equation}
   a_{n,\nu} =
   \sqrt{\frac{2}{\pi}}\frac{1}{\sqrt{nh''(t_{n,\nu})}}\cos\left(\rho_{n,\nu}+\frac{\pi}{4}\right)
   + \widetilde{R}_{n,\nu}
   \label{eq:central_terms_result}
\end{equation}
where $\rho_{n,\nu} = n(h(t_{n,\nu}) - h'(t_{n,\nu})t_{n,\nu})$ and 
\begin{equation}
   \abs{\widetilde{R}_{n,\nu}} \le \frac{C}{nh''(t_{n,\nu})\delta_n} +
   n\omega(\delta_n)\delta_n^3.
   \label{eq:central_terms_remainder}
\end{equation}

We note, in passing, that $\rho_{n,\nu} = -nh^{*}\left( \frac{\nu}{n} \right)$, where $h^{*}$ is the Legendre transform of the function $h$, i.e.~$h^*(x) = x\psi(x) - h(\psi(x))$, $\psi$ being the inverse function of the derivative of $h$.  

We write
\begin{equation}
   b_{n,\nu} = B_{n,\nu} + R_{n,\nu}^{(1)}
   \label{eq:central_terms_split}
\end{equation}
where
\begin{equation}
   B_{n,\nu} = \int_{t_{n,\nu}-\delta_n}^{t_{n,\nu}+\delta_n} \exp\left( i\left( nh(t) - \nu t \right) \right) dt
   \label{eq:central_bnnu}
\end{equation}
and
\begin{equation}
   R_{n,\nu}^{(1)} = 
   \int_{0}^{t_{n,\nu}-\delta_n} \exp\left( i\left( nh(t) - \nu t \right) \right) dt
   +\int_{t_{n,\nu}+\delta_n}^\pi \exp\left( i\left( nh(t) - \nu t \right)
   \right) dt.
   \label{eq:central_Rnnu}
\end{equation}

Let $\varphi(t) = nh(t) - \nu t$, then since $h'$ is increasing we have for $t
\in \left[0, t_{n,\nu}-\delta_n\right]$ that $\varphi'(t) = nh'(t) - \nu =
n(h'(t) - h'(t_{n,\nu})) \le n\left( h'(t_{n,\nu} - \delta_n) - h'(t_{n,\nu})
\right) = -n\delta_nh''(\xi)$ for some $\xi$ between $t_{n,\nu}-\delta_n$ and
$t_{n,\nu}$. So by part (b) of Lemma~\ref{lem:lem1}, and then by case
\ref{case:useful_ineq_c3} of inequality \eqref{eq:useful_hdoubleprime_inequality} we obtain
\[
\abs{\int_0^{t_{n,\nu}-\delta_n} \exp\left( i(nh(t) - \nu t) \right) dt} \le \frac{2}{nh''(\xi)\delta_n} \le \frac{C}{n\delta_nh''(t_{n,\nu})}.
\]
Same type of estimate holds for $\int_{t_{n,\nu} + \delta_n}^\pi \exp\left( i(nh(t) - \nu t) \right) dt$, so we obtain
\begin{equation}
   \abs{R_{n,\nu}^{(1)}} \le \frac{C}{nh''(t_{n,\nu})\delta_n}.
   \label{eq:rnnu_estim}
\end{equation}

Now we consider the integral $B_{n,\nu}$ defined in \eqref{eq:central_bnnu}.  Writing $h(t) = h(t_{n,\nu}) + h'(t_{n,\nu})(t - t_{n,\nu}) + \frac{h''(\xi)}{2}(t-t_{n,\nu})^2$ for $\abs{t-t_{n,\nu}} < \delta_n$ and some $\xi$ between $t$ and $t_{n,\nu}$, and using $h'(t_{n,\nu}) = \frac{\nu}{n}$, we can easily check that
\begin{equation}
   \exp\left( i(nh(t) - \nu t) \right) = e^{i\rho_{n,\nu}}e^{inh''(t_{n,\nu})(t-t_{n,\nu})^2/2} + R_{n,\nu}^{(2)}(t)
   \label{eq:central_taylor_exp}
\end{equation}
where
\[
\begin{split}
\abs{R_{n,\nu}^{(2)}(t)} & = \abs{e^{i\frac{n}{2}\left( h''(\xi) - h''(t_{n,\nu})  \right)(t-t_{n,\nu})^2} - 1} \\
& \le \frac{n}{2}\abs{h''(\xi) - h''(t_{n,\nu})}(t-t_{n,\nu})^2 \\
& \le \frac{n}{2}\omega(\delta_n)\delta_n^2.
\end{split}
\]
From the last estimate and \eqref{eq:central_taylor_exp} we obtain that
\begin{equation}
   B_{n,\nu} = \widetilde{B_{n,\nu}} + R_{n,\nu}^{(2)}
   \label{eq:split_of_Bnnu}
\end{equation}
where
\begin{equation}
   \widetilde{B_{n,\nu}} = e^{i\rho_{n,\nu}} \int_{t_{n,\nu}-\delta_n}^{t_{n,\nu}+\delta_n} \exp\left( i\frac{h''(t_{n,\nu})}{2}n(t-t_{n,\nu})^2 \right) dt
   \label{eq:tilde_bnnu_def}
\end{equation}
and
\begin{equation}
   \abs{R_{n,\nu}^{(2)}} = \abs{\int_{t_{n,\nu}-\delta_n}^{t_{n,\nu}+\delta_n} R_{n,\nu}^{(2)}(t) dt} \le n\omega(\delta_n)\delta_n^3.
   \label{eq:estim_R2nnu}
\end{equation}
A change of variable $u = \sqrt{\frac{nh''(t_{n,\nu})}{2}}(t-t_{n,\nu})$ transforms \eqref{eq:tilde_bnnu_def} into
\begin{equation}
   \begin{split}
      \widetilde{B_{n,\nu}} & = \frac{\sqrt{2}e^{i\rho_{n,\nu}}}{\sqrt{nh''(t_{n,\nu})}}\int_{-\delta_n\sqrt{\frac{nh''(t_{n,\nu})}{2}}}^{\delta_n\sqrt{\frac{nh''(t_{n,\nu})}{2}}} e^{iu^2} du \\
      & = \frac{2\sqrt{2}e^{i\rho_{n,\nu}}}{\sqrt{nh''(t_{n,\nu})}}\left(\int_0^\infty e^{iu^2} du - \int_{\delta_n\sqrt{\frac{nh''(t_{n,\nu})}{2}}}^{\infty} e^{iu^2} du\right).
   \end{split}
   \label{eq:widetilde_Bnnu_substitution}
\end{equation}
Since $\int_0^\infty e^{iu^2} du = \frac{\sqrt{\pi}}{2}e^{i\pi/4}$ and 
\[
\abs{\int_x^\infty e^{iu^2} du} = \abs{\int_{x^2}^\infty \frac{e^{it}}{2\sqrt{t}} dt} \le \frac{1}{x}, 
\]
we get from \eqref{eq:widetilde_Bnnu_substitution} that
\begin{equation}
   \widetilde{B_{n,\nu}} = \frac{\sqrt{2\pi}e^{i\left(\rho_{n,\nu}+\frac{\pi}{4}\right)}}{\sqrt{nh''(t_{n,\nu})}} + R_{n,\nu}^{(3)}
   \label{eq:widetilde_Bnnu_almost_done}
\end{equation}
where
\begin{equation}
   \abs{R_{n,\nu}^{(3)}} \le \frac{4}{nh''(t_{n,\nu})}\frac{1}{\delta_n}.
   \label{eq:Rnnu3_estim}
\end{equation}
From \eqref{eq:central_terms_split}, \eqref{eq:split_of_Bnnu} and \eqref{eq:widetilde_Bnnu_almost_done} we obtain
\begin{equation}
   b_{n,\nu} = \frac{\sqrt{2\pi}}{\sqrt{nh''(t_{n,\nu})}}e^{i\left( \rho_{n,\nu}+\frac{\pi}{4} \right)} + \widetilde{\widetilde{R_{n,\nu}}},
   \label{eq:bnnu_final}
\end{equation}
where $\widetilde{\widetilde{R_{n,\nu}}} = R_{n,\nu}^{(1)} + R_{n,\nu}^{(2)} + R_{n,\nu}^{(3)}$, and, by \eqref{eq:rnnu_estim}, \eqref{eq:estim_R2nnu} and \eqref{eq:Rnnu3_estim} we get
\begin{equation}
   \abs{\widetilde{\widetilde{R_{n,\nu}}}} \le \frac{C}{nh''(t_{n,\nu})}\frac{1}{\delta_n} + n\omega(\delta_n)\delta_n^3.
   \label{eq:Rtotal_estim}
\end{equation}
Since $h$ is odd, 
\[
a_{n,\nu} = \frac{2 \Rp\left( b_{n,\nu}\right)}{2\pi}
\]
and, by \eqref{eq:bnnu_final}
\begin{equation}
   a_{n,\nu} = \sqrt{\frac{2}{\pi}}\frac{1}{\sqrt{nh''(t_{n,\nu})}}\cos\left( \rho_{n,\nu} + \frac{\pi}{4} \right) + \widetilde{R_{n,\nu}},
   \label{eq:annu_final}
\end{equation}
where $\abs{\widetilde{R_{n,\nu}}} \le \abs{\widetilde{\widetilde{R_{n,\nu}}}}$.  By summing over the central terms, we obtain
\begin{equation}
   \sum_{n\alpha_n \le \nu \le n\beta_n} \abs{a_{n,\nu}} = \sqrt{\frac{2}{\pi}} \sum_{n\alpha_n \le \nu \le n\beta_n} \frac{1}{\sqrt{nh''(t_{n,\nu})}}\abs{\cos\left( \rho_{n,\nu} + \frac{\pi}{4} \right)} + R_n
   \label{eq:central_terms_sum}
\end{equation}
where
\begin{equation}
   \abs{R_n} \le \sum_{n\alpha_n \le \nu \le n\beta_n} \abs{\widetilde{\widetilde{R_{n,\nu}}}} \le \frac{C}{\delta_n}\sum_{n\alpha_n \le \nu \le n\beta_n} \frac{1}{nh''(t_{n,\nu})} + Cn^2\omega(\delta_n)\delta_n^3
   \label{eq:central_sum_remainder_estim}
\end{equation}
(Here we have used the fact that the number of central terms is $(\beta_n - \alpha_n)n < Cn$.)

We observe that 
\[
   \sum_{n\alpha_n \le \nu \le n\beta_n} \frac{1}{nh''(t_{n,\nu})} = \sum_{n\alpha_n \le \nu \le n\beta_n} \frac{1}{nh''\left( \psi\left( \frac{\nu}{n} \right) \right)}
\]
is bounded, since it is positive and smaller than 
\begin{multline*}
\sum_{n\alpha + 1 \le \nu \le n\beta-1} \frac{1}{nh''\left( \psi\left( \frac{\nu}{n} \right) \right)} = \\ 
\sum_{n\alpha + 1 \le \nu \le nh'\left( \frac{\pi}{2} \right)} \frac{1}{nh''\left( \psi\left( \frac{\nu}{n} \right) \right)} + \sum_{nh'\left( \frac{\pi}{2} \right) < \nu \le n\beta-1} \frac{1}{nh''\left( \psi\left( \frac{\nu}{n} \right) \right)} 
\end{multline*}
and the last two summands converge to a finite limit. (Here we appealed to
Lemma~\ref{lem:part_riem_sum}, the conditions in that lemma being
satisfied by observations made in sections~\ref{prop:first}
and~\ref{prop:second})
Therefore \eqref{eq:central_sum_remainder_estim} gives
\[
\abs{R_n} \le C\left( \frac{1}{\delta_n} + n^2\omega(\delta_n)\delta_n^3 \right).
\]
Choosing $\delta_n$ so that $\frac{1}{\delta_n} = n^2\omega(\delta_n)\delta_n^3$, i.e.~choosing $\Phi_n$ so that
\begin{equation}
   \omega\left( \frac{\Phi_n}{\sqrt{n}} \right) \Phi_n^4 = 1,
   \label{eq:choice_of_Phi_n}
\end{equation}
the past estimate simplifies to 
\begin{equation}
   \abs{R_n} \le \frac{C}{\delta_n} = C\frac{\sqrt{n}}{\Phi_n}. 
   \label{eq:final_central_terms_estimate}
\end{equation}

\begin{remark}
   Since it is a property of a modulus of continuity that there exists $C > 0$
   such that $\omega(x) > Cx$, we obtain from \eqref{eq:choice_of_Phi_n} that
   $Phi_n = O(n^{1/10})$, $n \rightarrow \infty$.
\end{remark}

\section{Final Step}
Taking into account contributions of external \eqref{eq:sum_external} and periphery
\eqref{eq:periphery_sum} terms, together with \eqref{eq:central_terms_sum} and
\eqref{eq:final_central_terms_estimate} and the last remark, we obtain
\begin{equation}
   \frac{1}{\sqrt{n}}\sum_{\nu=-\infty}^\infty \abs{a_{n,\nu}} = \sqrt{\frac{2}{\pi}} \sum_{\nu=n\alpha_n}^{n\beta_n} \frac{1}{n}\frac{1}{\sqrt{h''(t_{n,\nu})}}\abs{\cos\left( \rho_{n,\nu}+\frac{\pi}{4} \right)} + \mathcal{O}\left( \frac{1}{\Phi_n} \right), n \rightarrow \infty.
   \label{eq:all_terms_sum}
\end{equation}
We shall prove that 
\begin{equation}
   \mathcal{D}_n = \sum_{\nu=n\alpha_n}^{n\beta_n} \frac{1}{n}
   \frac{1}{\sqrt{h''(t_{n,\nu})}}\abs{\cos\left(
   \rho_{n,\nu}+\frac{\pi}{4} \right)} - \frac{2}{\pi}\int_\alpha^\beta
   \frac{dt}{\sqrt{h''(\psi(t))}} \rightarrow 0, n \rightarrow \infty.
   \label{eq:final_steps_to_prove}
\end{equation}
Using the simplification of the last integral given in
\eqref{eq:int_1_over_sqrt_equals}, i.e.~making change of variable $x =
\psi(t)$, we obtain from \eqref{eq:all_terms_sum} and
\eqref{eq:final_steps_to_prove} the conclusion of Theorem~\ref{thm:main}

To prove \eqref{eq:final_steps_to_prove}, for each sufficiently small
$\varepsilon > 0$ we write
\[
\mathcal{D}_n = A_{n,\varepsilon} + B_{n,\varepsilon} + C_{n,\varepsilon} -
E_{\varepsilon}
\]
where
\[
A_{n,\varepsilon} = 
\sum_{\nu=n\alpha_n}^{nh'(\varepsilon)} \frac{1}{n}
\frac{1}{\sqrt{h''(t_{n,\nu})}}\abs{\cos\left(
\rho_{n,\nu}+\frac{\pi}{4} \right)},
\]
\[
C_{n,\varepsilon} = 
\sum_{\nu=nh'(\pi-\varepsilon)}^{n\beta_n} \frac{1}{n}
\frac{1}{\sqrt{h''(t_{n,\nu})}}\abs{\cos\left(
\rho_{n,\nu}+\frac{\pi}{4} \right)},
\]
\[
E_{\varepsilon} = 
\frac{2}{\pi}\int_{h'(0)}^{h'(\varepsilon)}
\frac{dt}{\sqrt{h''(\psi(t))}}
+ \frac{2}{\pi}\int_{h'(\pi-\varepsilon)}^{h'(\pi)}
\frac{dt}{\sqrt{h''(\psi(t))}},
\]
and
\[
B_{n,\varepsilon} = \sum_{\nu=nh'(\varepsilon)}^{nh'(\pi-\varepsilon)} \frac{1}{n}
\frac{1}{\sqrt{h''(t_{n,\nu})}}\abs{\cos\left(
\rho_{n,\nu}+\frac{\pi}{4} \right)} -
\frac{2}{\pi}\int_{h'(\varepsilon)}^{h'(\pi-\varepsilon)}
\frac{dt}{\sqrt{h''(\psi(t))}}.
\]
Since all the summands in $A_{n,\varepsilon}$ are non-negative, we get 
\[
0 \le A_{n,\varepsilon} \le \sum_{n\alpha+1 \le \nu \le nh'(\varepsilon)}
\frac{1}{n}\frac{1}{\sqrt{h''(t_{n,\nu})}}\cdot 1
\]
and, by \eqref{eq:par_3_estim1_for_sqrt} it follows that
$\abs{A_{n,\varepsilon}}\le C\varepsilon$, and, in the same way,
$\abs{C_{n,\varepsilon}}\le C\varepsilon$. 

We obtain similarly from \eqref{eq:int_1_over_sqrt_estim} that $0 \le
E_\varepsilon \le 2\varepsilon$. It follows then that
\begin{equation}
   \abs{\mathcal{D}_n} \le \abs{B_{n,\varepsilon}} + 2(C+1)\varepsilon.
   \label{eq:reduction_to_Bn}
\end{equation}

Now we let $\varphi(x) = \frac{1}{\pi} \left( h(\psi(x)) - x\psi(x)
\right)$ ($=-\frac{1}{\pi}h^*(x)$ where $h^*$ denotes the Legendre transform
of $h$), and observe that $\varphi'(x) = -\frac{1}{\pi}\psi(x)$, and
$\varphi''(x) = -\frac{1}{\pi} \psi'(x) =
-\frac{1}{\pi}\frac{1}{h''(\psi(x))}$.
Let $J_\varepsilon$ denote the interval $\left[ h'(\varepsilon),
h'(\pi-\varepsilon) \right]$.  For $x\in J_\varepsilon$ we have
$\psi(x)\in\left[ \varepsilon,\pi-\varepsilon \right]$, and it follows that
$\varphi$ is twice continuously differentiable on $J_\varepsilon$, and that its
second derivative is bounded away from zero on $J_\varepsilon$.  Therefore, by
Lemma~\ref{lem:equidist}, the array of fractional parts of $n\varphi\left(
\frac{\nu}{n} \right) = \frac{\rho_{n,\nu}}{\pi}$, $n=1,2,3,\dots$, $\frac{\nu}{n} \in J_\varepsilon$ is
equidistributed (in the sense of Lemma~\ref{lem:snk_equiv}). 

Now we shall apply Lemma~\ref{lem:snk_equiv} to the interval $J_\varepsilon$, the
array of fractional parts of $\frac{\rho_{n,\nu}}{\pi}$, and the functions
$f(x) = \frac{1}{\sqrt{h''(\psi(x))}}$ and $g(x) = \abs{\cos\pi\left(x +
\frac{1}{4} \right)}$, which are Riemann integrable on $J_\varepsilon$ and
$[0,1]$, respectively, $g$ being periodic with period $1$. Since
\[
\frac{1}{n} \sum_{\left\{ \nu|\frac{\nu}{n} \in J_\varepsilon \right\}}
f\left( \frac{\nu}{n} \right) g\left( \frac{\rho_{n,\nu}}{\pi} \right) =
\frac{1}{n} \sum_{\left\{ \nu|\frac{\nu}{n} \in J_\varepsilon \right\}}
\frac{1}{\sqrt{h''\left( \psi\left( \frac{\nu}{n} \right) \right)}}
\abs{\cos\left( \rho_{n,\nu}+\frac{\pi}{4} \right)},
\]
we obtain by Lemma~\ref{lem:snk_equiv}(c) that the last expression tends to 
\[
\int_{J\varepsilon} \frac{1}{\sqrt{h''(\psi(x))}} dx \cdot
\int_0^1\abs{\cos\pi\left( x + \frac{1}{4} \right)}dx =
\frac{2}{\pi}\int_{J_\varepsilon} \frac{1}{\sqrt{h''(\psi(x))}} dx.
\]
Therefore $B_{n,\varepsilon} \rightarrow 0$ as $n\rightarrow\infty$ for any
$\varepsilon > 0$, so from \eqref{eq:reduction_to_Bn} we get $\mathcal{D}_n
\rightarrow 0$ as $n \rightarrow\infty$.  This proves
\eqref{eq:final_steps_to_prove} and therefore Theorem~\ref{thm:main}.

\section{Proof of Theorem~\ref{thm:stronger}}

\begin{proof}
   Let 
   \[
   S(x) = S(x,h) = \frac{1}{\sqrt{x}} \sum_{\nu=-\infty}^\infty \abs{a_{x,\nu}} = \frac{1}{\sqrt{x}}\norm{e^{ixh(t)}}
   \]
   where
   \[
   a_{x,\nu} = \frac{1}{2\pi}\int_{-\pi}^\pi \exp\left( i(xh(t) - \nu t \right)dt,
   \]
   and let
   \[
   L = L(h) = \left( \frac{2}{\pi} \right)^{\frac{3}{2}} \int_0^\pi \sqrt{\abs{h''(t)}}dt.
   \]
   We shall prove that
   \begin{enumerate}[label=(\arabic*)]
      \item $S(x)$ is continuous on $(1,\infty)$, and \label{case:cond1}
      \item $S(n\alpha) \rightarrow L$ as $n \rightarrow \infty$ for any $\alpha$ real, $\alpha\neq 0$. \label{case:cond2}
   \end{enumerate}
   It is a well known fact, obtained by a category argument, that \ref{case:cond1} and \ref{case:cond2} imply that $S(x) \rightarrow L$ as the real variable $x$ tends to $\infty$. Therefore, to prove Theorem~\ref{thm:stronger}, it is sufficient to prove \ref{case:cond1} and \ref{case:cond2}. 

   Applying Lemma~\ref{lem:lem1} part \ref{lem:lem1_part_a} to the case $\varphi(t) = xh(t)-\nu t$ on $\left[ -\pi,\pi \right]$, we see that for any $A>0$ there exists $C>0$ such that if $1 \le \abs{x}\le A$, $\abs{\nu}\ge 2A\max\abs{h'}$, then $\abs{a_{x,\nu}} \le \frac{C}{\nu^2}$. By the Weierstrass M-test the series $\sum \abs{a_{x,\nu}}$ is uniformly convergent on $\left[ -A,A \right]$, and since $a_{x,\nu}$ are continuous function of $x$ for all $\nu$, $S(x,h)$ is continuous on $(1,\infty)$.  

   If $h$ is $2\pi$-periodic, so is $\alpha h$, and $\alpha h$ satisfies the conditions of Theorem~\ref{thm:main} for all $\alpha>0$. So, by Theorem~\ref{thm:main}, $S(n,\alpha h) \rightarrow L(\alpha h)$ and we obtain $S(n\alpha, h) = \frac{1}{\sqrt{\alpha}}S(n,\alpha h) \rightarrow \frac{1}{\sqrt{\alpha}} L(\alpha h) = L(h)$, which proves \ref{case:cond2}.
\end{proof}

\section{Corollaries}
\begin{cor} (Oral communication of George Stey)
   Let $J_\nu$ denotes the Bessel function of integer order $\nu$.  Then
   \[
   \frac{1}{\sqrt{x}}\sum_{\nu=-\infty}^\infty \abs{J_\nu(x)} \rightarrow
   \frac{16}{\Gamma^2\left( \frac{1}{4} \right)}
   \]
   as $x$ tends to $\infty$.
\end{cor}
\begin{proof}
   From the well known fact
   \[e^{ix\sin(t)} = \sum_{\nu=-\infty}^\infty J_\nu(x)e^{i\nu t}\]
   we obtain that 
   \[
   \sum_{\nu=-\infty}^\infty\abs{J_\nu(x)} = \norm{e^{ix\sin(t)}}
   \]
   Since $h(t) = \sin(t)$ satisfies the conditions of Theorem~\ref{thm:stronger}, and since
   \[
   \int_0^\pi \sqrt{\abs{h''(t)}}dt = \int_0^\pi\sqrt{\sin(t)}dt = B\left(
   \frac{3}{4}, \frac{1}{2} \right) = \frac{4\pi\sqrt{2\pi}}{\Gamma^2\left(
   \frac{1}{4} \right)}
   \]
   we obtain from the Theorem~\ref{thm:stronger} that 
   \[
   \frac{1}{\sqrt{x}}\sum_{\nu=-\infty}^\infty \abs{J_\nu(x)} \rightarrow
   \frac{16}{\Gamma^2\left( \frac{1}{4} \right)} \text{, $x \rightarrow\infty$.}
   \]
\end{proof}

\begin{cor}\label{col:product}
   Let
   \begin{equation}
      f(t) = \prod_{j = 1}^{J} \frac{e^{it} - \alpha_j}{1 - \alpha_j e^{it}}
   \end{equation}
   where $\alpha_j \in (0,1)$ for all $j = 1,\dots,J$.

   Then 
   \begin{equation}\label{eq_example}
      \frac{\norm{f^n}}{\sqrt{n}} = \left(\frac{2}{\pi}\right)^{\frac{3}{2}}
      \int_{-1}^{1} \sqrt{\sum_{j=1}^{J} \frac{2 \alpha_j \left(1 -
      \alpha_j^2\right)}{\left(1 + \alpha_j^2 - 2\alpha_j u\right)^2}}
      \left(1 - u^2\right)^{-\frac{1}{4}} du +
      o(1) \text{as $n \rightarrow\infty$}
   \end{equation}
\end{cor}
\begin{proof}
   It is obvious that $f$ is continuous $2\pi$-periodic infinitely differentiable
   function, and it is easy to see that $f(t) = e^{-ih(t)}$, where
   \begin{equation}
      h(t) = \sum_{j=1}^{J} i\log\left(\frac{e^{it} - \alpha_j}{1 -
      \alpha_je^{it}}\right).
   \end{equation}
   The first three derivatives of $h$ are
   \begin{subequations}
      \begin{align}
	 h'(t)&  = -\sum_{j=1}^{J} \frac{1 - \alpha_j^2}{1 + \alpha_j^2 - 2\alpha_j\cos t}
	 \label{eq_h_firstder}\\
	 h''(t)& = \sum_{j=1}^{J} \frac{2\alpha_j\left(1 - \alpha_j^2\right)\sin
	 t}{\left(1 + \alpha_j^2 - 2\alpha_j\cos t\right)^2}
	 \label{eq_h_secder}\\
	 h'''(t)& = \sum_{j=1}^{J} \frac{2\alpha_j\left(1 -
	 \alpha_j^2\right)\left(-3 \alpha_j + \cos(t) + \alpha_j^2 \cos(t) + \alpha_j
	 \cos(2t)\right)}{\left(1 + \alpha_j^2 - 2 \alpha_j \cos(t)\right)^3}
      \end{align}
   \end{subequations}
   As each of the terms
   of the sum in~\eqref{eq_h_secder} is positive,  it is obvious that $h''(t) >
   0$ for $t \in (0,\pi)$.  Finally,
   \begin{align}
      h'''(0)& = \sum_{j=1}^{J} \frac{2\alpha_j\left(1 +
      \alpha_j\right)}{\left(1 - \alpha_j\right)^3} \neq 0\\
      h'''(\pi)& = \sum_{j=1}^{J} \frac{2\alpha_j\left(1 -
      \alpha_j\right)}{\left(1 + \alpha_j\right)^3} \neq 0.
   \end{align}

   Therefore the function $h$ satisfies all the conditions of the Theorem~\ref{thm:main}, and
   since $\norm{f^n} = \norm{e^{-inh(t)}} = \sum_{\nu=-\infty}^\infty
   \abs{a_{n,\nu}}$, where $f^n(t) = \sum_{\nu=-\infty}^\infty
   a_{n,\nu}e^{i\nu t}$, we have
   \begin{equation}
      \frac{1}{\sqrt{n}}\sum_{\nu = -\infty}^{\infty} \abs{a_{n,\nu}} =
      \left(\frac{2}{\pi}\right)^{\frac{3}{2}} \int_{0}^{\pi}
      \sqrt{\sum_{j=1}^{J} \frac{2\alpha_j\left(1 - \alpha_j^2\right)\sin
      t}{\left(1 + \alpha_j^2 - 2\alpha_j\cos t\right)^2}} dt +
      o(1) \text{ as $n\rightarrow\infty$}
   \end{equation}
   Substituting $u = \cos(t)$ and simplifying will result in~\eqref{eq_example}.
\end{proof}
\begin{cor} (Girard's formula)
   If 
   \[g(z) = \frac{z-\alpha}{1-\overline{\alpha}z}, 0 < \abs{\alpha} < 1\]
   then
   \begin{equation}\label{eq:colgir}
      \frac{\norm{g^n\left( e^{it} \right)}}{\sqrt{n}} = 16 \sqrt{2}
      \left(\Gamma\left(\frac{1}{4}\right)\right)^{-2} \frac{\sqrt{\abs{\alpha}}}{1 +
      \abs{\alpha}}\; F\left(\frac{1}{2}, \frac{3}{4}; \frac{3}{2};
      \frac{4\abs{\alpha}}{(1+\abs{\alpha})^2}\right) + 
      o(1) 
   \end{equation}
   as $n\rightarrow\infty$.
\end{cor}

\begin{proof}
Without loss of generality we may assume that $\alpha$ is real and positive.
Setting in Corollary\ref{col:product} $J=1$, $f(t) = g\left( e^{it} \right)$
and $\alpha_1 = \alpha$, we obtain
\begin{equation}\label{eq:colprodresult}
   \frac{\norm{f^n}}{\sqrt{n}} = \left(\frac{2}{\pi}\right)^{\frac{3}{2}}
   \int_{-1}^{1} \sqrt{\frac{2 \alpha \left(1 -
   \alpha^2\right)}{\left(1 + \alpha^2 - 2\alpha u\right)^2}}
   \left(1 - u^2\right)^{-\frac{1}{4}} du +
   o(1) \text{ as $n\rightarrow\infty$.}
\end{equation}

Substituting $s = \frac{1+u}{2}$ the integral will become
\begin{equation}\label{eq_int_to_hyper}
   2\frac{\sqrt{\alpha \left(1-\alpha^2\right)}}{\left(1+\alpha\right)^2}
   \int_{0}^{1} \left(1 - \frac{4\alpha}{(1+\alpha)^2} s\right)^{-1} 
   \left(1 - s\right)^{-\frac{1}{4}}s^{-\frac{1}{4}} \mathrm{d}s 
\end{equation}
which can be written as
\begin{equation}
   8\pi^{\frac{3}{2}}\left(\Gamma\left(\frac{1}{4}\right)\right)^{-2} \frac{\sqrt{\alpha}}{1 +
   \alpha}\; F\left(\frac{1}{2}, \frac{3}{4}; \frac{3}{2};
   \frac{4\alpha}{(1+\alpha)^2}\right)
\end{equation}
(see for example theorems 16  and 21 in~\cite{rainville}) 
Substituting this for the integral in~\eqref{eq:colprodresult} will give
us~\eqref{eq:colgir}.
\end{proof}

\section{Two Open Questions}

Theorem~\ref{thm:main} and Corollary~\ref{col:product} make us wonder to what extent the convergence of the sequence $\frac{1}{\sqrt{n}}\norm{\exp(inh(t))}$ is an exceptional phenomenon, and the following two questions arise:
\begin{enumerate}
   \item If the set of conditions in Theorem~\ref{thm:main} is weakened so that the condition \ref{cond2}: $h$ is odd is replaced by the condition (2'): $h''$ has no zeroes in $(-\pi,0)$, is it still true that $\frac{1}{\sqrt{n}}\norm{\exp(inh(t))}$ converges?
   \item If $B(z)$ is a finite Blaschke product, is it true that $\frac{1}{\sqrt{n}}\norm{B^n\left( e^{it} \right)}$ converges?
\end{enumerate}

\bibliographystyle{plain}

\end{document}